\documentclass[11pt]{amsart}
\pdfoutput=1

\usepackage{amsmath,amsthm, amssymb, amsfonts, epsfig}
\usepackage[mathscr]{eucal}
\usepackage{graphics}
\usepackage{color}
\usepackage{subfigure}

\newtheorem{theorem}{Theorem}[section]

\newtheorem{lemma}[theorem]{Lemma}

\theoremstyle{definition}
\newtheorem{remark}[theorem]{Remark}

\newcommand{\thmref}[1]{Theorem~\ref{#1}}

\newcommand{\lemref}[1]{Lemma~\ref{#1}}

\newcommand{\mr}[1]{\mathrm{#1}}

\usepackage[colorlinks=false, pdftitle={Higher-dimensional linking integrals}, pdfauthor={Clayton Shonkwiler, David Shea Vela-Vick}, pdfsubject={Linking integrals}, pdfkeywords={Gauss linking integral, linking number, 57Q45, 57M25, 53C20}]{hyperref}
\newtheorem{fact}[theorem]{Fact}
\begin{document}

\thispagestyle{empty}

\title{Higher-Dimensional Linking Integrals}
\author{Clayton Shonkwiler \and David Shea Vela-Vick}
\address{Department of Mathematics \\ University of Pennsylvania}
\email{shonkwil@math.upenn.edu}
\urladdr{http://www.math.upenn.edu/~shonkwil}
\address{Department of Mathematics \\ University of Pennsylvania}
\email{dvick@math.upenn.edu}
\urladdr{http://www.math.upenn.edu/~dvick}

\date{January 25, 2008}
\maketitle

\begin{abstract}
We derive an integral formula for the linking number of two submanifolds of the n-sphere $S^n$, of the product $S^n \times \mathbb{R}^m$, and of other manifolds which appear as ``nice'' hypersurfaces in Euclidean space.  The formulas are geometrically meaningful in that they are invariant under the action of the special orthogonal group on the ambient space.
\end{abstract}

\section{Introduction} 
	\label{sec:intro}
	In \cite{Gauss}, Carl Friedrich Gau\ss \ defined an integral whose purpose was to compute the linking number of two closed curves in Euclidean 3-space.  Gau\ss's linking integral follows from at least two elementary arguments: one involves thinking of the curves as wires, running a current through one and applying Amp\`ere's Law, while the other is a straightforward degree-of-map argument (see \cite{Epple}).  The latter extends easily to all Euclidean spaces (as we will see in Section \ref{sec:rN}), but cannot be directly adapted to spheres.  As such, when Dennis DeTurck and Herman Gluck set out to derive integral formulas for the linking number of closed curves in the 3-sphere and in hyperbolic 3-space \cite{DG1}, they adapted the first argument, developing steady-state versions of electrodynamics in those spaces.  
	
	The integrals defined by Gau\ss \ in the case of $\mathbb{R}^3$ and by DeTurck and Gluck in the case of $S^3$ and $H^3$ apply to the simplest of the eight 3-dimensional model geometries and motivate the search for a linking integral in the next simplest geometry, $S^2 \times \mathbb{R}$.  An explicit formula for this case is given in Section \ref{sec:examples}, but the techniques developed to obtain this formula apply to a much broader class of manifolds, which we call ``visible hypersurfaces''.

	Call a smooth hypersurface $M^n \subset \mathbb{R}^{n+1}$ {\it visible from the point $p$} if each ray from $p$ either misses $M^n$ completely, or else meets it just once transversally.  We will always arrange things so that $p$ is the origin in $\mathbb{R}^{n+1}$, and simply call our hypersurfaces {\it visible}.  The principal examples are $S^n \subset \mathbb{R}^{n+1}$ and $S^n \times \mathbb{R}^m \subset \mathbb{R}^{n+m+1}$.  Every closed manifold is homotopically equivalent to a visible hypersurface in some Euclidean space.  Our main theorem gives explicit linking integrals in all visible hypersurfaces:

\begin{theorem}\label{thm:main}
	Let $K^k$ and $L^\ell$ be disjoint, closed, oriented, null-homologous submanifolds of a visible hypersurface $M^n \subset \mathbb{R}^{n+1}$ such that $k + \ell = n-1$.  Then their linking number is given by
	\begin{equation}\label{eqn:snxrmIntegral}
		Lk(K, L) = \frac{1}{\mathrm{vol}\, S^{n}} \int_{K \times L} \frac{\Omega_{k, \ell}(\alpha)}{|x|^{k+1} |y|^{\ell+1} \sin^{n}\alpha} \, [x, dx, y, dy]
	\end{equation}
	where
	\[
		\Omega_{k, \ell}(\alpha) = \int_{\theta=\alpha}^\pi \sin^k (\theta - \alpha) \sin^\ell \theta \, d\theta.
	\]
	Here $\alpha(x, y)$ is the angle between $x \in K$ and $y \in L$, thought of as vectors in $\mathbb{R}^{n+1}$, and the notation $[x,dx, y, dy]$ is as defined in Section \ref{sec:sn}.
\end{theorem}

In $\mathbb{R}^3$, or equivalently $S^3$, the linking number of $K$ with $L$ is defined to be the oriented intersection number of $K$ with a chain $\overline{L}$ bounded by $L$.  Extend this definition to higher dimensions by defining the linking number of two submanifolds $K^k$ and $L^\ell$ to be the oriented intersection number $K \cdot \overline{L}$ of $K$ with a chain $\overline{L}$ bounded by $L$.  This intersection number only makes sense when $k$ and $\ell$ are as in \thmref{thm:main}.

\begin{remark}\label{rem:affine}
	The integrand in \eqref{eqn:snxrmIntegral} is $SO(n+1)$-invariant, meaning that it will be the same for $K$ and $L$ in $M^n \subset \mathbb{R}^{n+1}$ as it is for $h(K)$ and $h(L)$ in $h(M)$ for all $h \in SO(n+1)$, whether or not $M^n$ is invariant under $h$.
\end{remark}

As an immediate corollary to \thmref{thm:main}, we get an integral formula for the linking number of submanifolds of $S^n$ which agrees with the formulas obtained independently and by quite different methods by Kuperberg \cite{Kuperberg} and DeTurck and Gluck \cite{DG2}:

\begin{theorem}\label{thm:sn}
	Let $K^k$ and $L^\ell$ be disjoint, closed, oriented submanifolds of the round $n$-sphere $S^n$ with $k + \ell = n-1$.  Then
	\[
		Lk(K, L) = \frac{1}{\mathrm{vol}\,  S^n} \int_{K \times L} \frac{\Omega_{k, \ell}(\alpha)}{\sin^n \alpha} \, [x, dx, y, dy]
	\]
	where
	\[
		\Omega_{k, \ell} (\alpha)= \int_{\theta = \alpha}^\pi \sin^k(\theta - \alpha) \sin^\ell \theta \, d\theta.
	\]
	Here $\alpha(x,y)$ is the distance in $S^n$ from $x \in K$ to $y \in L$.
\end{theorem}

\begin{remark}\label{rem:geometricallyMeaningful}
	The integrand in \thmref{thm:sn} is invariant under orientation-preserving isometries of $S^n$.
\end{remark}

To prove \thmref{thm:main}, we first adapt the degree-of-map proof of the Gau\ss \ linking integral to get an integral formula for the linking number of closed submanifolds of $\mathbb{R}^{N}$ for all $N$.  Then, for submanifolds $K$ and $L$ of the visible hypersurface $M^n \subset \mathbb{R}^{n+1}$, we associate to $K$ a family of closed singular submanifolds $CK_R \subset \mathbb{R}^{n+1}$ of one higher dimension for $R \in (1,\infty)$.  We apply the integral formula in $\mathbb{R}^{n+1}$ to the pair $(CK_R,L)$ and, after taking an appropriate limit, deduce the integral in \thmref{thm:main}.

\subsection*{Acknowledgments} 
\label{sub:acknowledgements}
	We owe a considerable debt of gratitude to Dennis DeTurck and Herman Gluck for introducing us to linking integrals and for many helpful and illuminating conversations.

\section{A linking integral for $\mathbb{R}^N$} 
\label{sec:rN}
	To obtain a linking integral formula for $\mathbb{R}^N$, we directly adapt the degree-of-map proof of the Gau\ss \ linking integral.
	
	Let $K^k, L^\ell \subset \mathbb{R}^N$ be disjoint, closed, oriented submanifolds such that $k + \ell = N-1$.  Let $x: \mathbb{R}^k \to K \subset \mathbb{R}^N$ and $y: \mathbb{R}^\ell \to L \subset \mathbb{R}^N$ be oriented local coordinates for $K$ and $L$, where ${\bf s} = (s_1, \ldots , s_k)$ and ${\bf t} = (t_1, \ldots , t_\ell)$ give the coordinates on $\mathbb{R}^k$ and $\mathbb{R}^\ell$, respectively.  Then, up to sign, we can express $Lk(K, L)$ as the degree of a map:
	\begin{lemma}\label{lem:degree}
		If $K$ and $L$ are as above and $f: K \times L \to S^{N-1}$ is given by
		\[
			f(x, y) = \frac{x-y}{|x-y|},
		\]
		then the degree of the map $f$ is equal to $(-1)^N$ times the linking number of $K$ and $L$.  In other words, if $\omega$ is a volume form on $S^{N-1}$, then 
		\[
			Lk(K, L) = \frac{(-1)^N}{\mathrm{vol}\,  S^{N-1}} \int_{K \times L} f^* \omega.
		\]
	\end{lemma} 
	The idea of the proof of \lemref{lem:degree} is the following: if we change $K$ by a homology in the complement of $L$ (or {\it vice versa}), the degree of $f$ remains unchanged.  Hence, we can replace $K$ and $L$ by pairwise meridional round $k$- and $\ell$-spheres (with multiplicity). In this special case, it is a straightforward exercise to check that the degree of $f$ is $(-1)^N$ times the linking number.  
		
	With \lemref{lem:degree} in hand, we turn to finding a linking integral for submanifolds of $\mathbb{R}^{N}$.  In the following theorem and throughout the rest of this paper we will use the notation $(v_1, \ldots , v_{N})$ to denote the $N\times N$ matrix with rows given by the vectors $v_1, \ldots , v_{N} \in \mathbb{R}^{N}$.
	\begin{theorem}\label{thm:rN}
		With $K$ and $L$ as above, 
		\begin{equation}\label{eqn:rn+1Integral}
			Lk(K, L) = \frac{(-1)^{k+1}}{\mathrm{vol}\, S^{N-1}}\int_{K \times L} \frac{1}{|x-y|^{N}} \, [x-y, dx, dy],
		\end{equation}
		where  
		\[
			[x-y, dx, dy] = \det \left(x-y, \frac{\partial x}{\partial s_1} , \cdots , \frac{\partial x}{\partial s_k}, \frac{\partial y}{\partial t_1} , \cdots , \frac{\partial y}{\partial t_\ell} \right) d{\bf s} \, d{\bf t}.
		\]
	\end{theorem}
	\begin{remark}\label{rem:tripleProduct}
		\thmref{thm:rN} is a direct generalization of the original Gau\ss \ linking integral,
		$$Lk(K,L) = \frac{1}{4\pi} \int_{K \times L} \frac{x-y}{|x-y|^3}  \cdot \left( \frac{dx}{ds} \times \frac{dy}{dt} \right) ds \, dt,$$
		since $\frac{1}{|x-y|^3}[x-y, dx, dy]$ is equal to the triple product in the above integrand.
	\end{remark}

	\begin{proof}\label{pf:lemrn+1}
		By \lemref{lem:degree}, it suffices to show that
		\begin{equation}\label{eqn:rn+1IntegralEquivalent}
			f^* \omega = \frac{(-1)^\ell}{|x-y|^{N}} \, [x-y, dx, dy]
		\end{equation}
		where $\omega$ is the volume form on $S^{N-1}$.
		
		For $S^{N-1}$ embedded in $\mathbb{R}^N$ in the usual way, we can give an explicit formula for the volume form $\omega$ on $S^{N-1}$: if $p \in S^{N-1}$ and $V_1, \ldots , V_{N-1} \in T_pS^{N-1}$, then $\omega_p(V_1, \ldots, V_{N-1})$ is equal to the volume of the parallelpiped spanned by the vectors $p, V_1, \ldots , V_{N-1}$ in $\mathbb{R}^N$:
		\[
			\omega_p(V_1, \ldots , V_{N-1}) = \det(p, V_1, \ldots , V_{N-1}).
		\]
		
		Therefore, 
		\begin{equation}\label{eqn:degreedet}
			f^*\omega\left(\frac{\partial}{\partial s_1}, \ldots , \frac{\partial}{\partial s_k}, \frac{\partial}{\partial t_1}, \ldots , \frac{\partial}{\partial t_\ell}\right)  = \det\left(f(x,y), \frac{\partial f}{\partial s_1}, \ldots, \frac{\partial f}{\partial s_k}, \frac{\partial f}{\partial t_1}, \ldots , \frac{\partial f}{\partial t_\ell}\right).
		\end{equation}
		Note that, if $f(x,y) = (f_1(x,y),\dots,f_N(x,y))$, then
		\begin{equation}\label{eqn:dfidsj}
			\frac{\partial f_i}{\partial s_\eta} = \frac{\partial}{\partial s_\eta}\left(\frac{x_i - y_i}{|x-y|}\right) = \frac{\frac{\partial x_i}{\partial s_\eta}}{|x-y|} + (x_i - y_i) \frac{\partial}{\partial s_\eta}\left(\frac{1}{|x-y|}\right)
		\end{equation}
		and
		\begin{equation}\label{eqn:dfidtj}
			\frac{\partial f_i}{\partial t_\nu} = \frac{\partial}{\partial t_\nu} \left(\frac{x_i - y_i}{|x-y|}\right) = \frac{-\frac{\partial y_i}{\partial t_\nu}}{|x-y|} + (x_i - y_i) \frac{\partial}{\partial t_\nu}\left(\frac{1}{|x-y|}\right).
		\end{equation}
		However, since, in \eqref{eqn:degreedet}, $f_i(x,y) = \frac{x_i - y_i}{|x-y|}$ is in the same column as both $\frac{\partial f_i}{\partial s_\eta}$ and $\frac{\partial f_i}{\partial t_\nu}$, the second terms in both \eqref{eqn:dfidsj} and \eqref{eqn:dfidtj} contribute nothing to the determinant.
		
		Thus, combining \eqref{eqn:degreedet}, \eqref{eqn:dfidsj} and \eqref{eqn:dfidtj}, we see that
		\[
			f^*\omega\left(\frac{\partial}{\partial s_1}, \ldots , \frac{\partial}{\partial s_k}, \frac{\partial}{\partial t_1}, \ldots , \frac{\partial}{\partial t_\ell}\right) = \frac{(-1)^\ell}{|x-y|^{N}}\det\left(x-y, \frac{\partial x}{\partial s_1}, \ldots , \frac{\partial x}{\partial s_k} ,  \frac{\partial y}{\partial t_1}, \ldots , \frac{\partial y}{\partial t_\ell} \right),
		\]
		which is just a restatement of \eqref{eqn:rn+1IntegralEquivalent}.
	\end{proof}


\section{A linking integral for visible hypersurfaces} 
\label{sec:sn}
In this section, we prove \thmref{thm:main} by making use of the linking integral in $\mathbb{R}^{n+1}$.  We assume $n >1$; in the special case that $n=1$, $K$ and $L$ are both zero-dimensional and \thmref{thm:main} can be verified directly.  To set notation, let $K^k$ and $L^\ell$ be disjoint, closed, oriented, null-homologous submanifolds of the visible hypersurface $M^n \subset \mathbb{R}^{n+1}$ such that $k + \ell = n-1$.  We use $x$, $y$, ${\bf s}$ and ${\bf t}$ as in Section \ref{sec:rN} and use $[x, dx, y, dy]$ to denote
\[
	\det\left(x, \frac{\partial x}{\partial s_1}, \ldots , \frac{\partial x}{\partial s_k}, y,  \frac{\partial y}{\partial t_1}, \ldots , \frac{\partial y}{\partial t_\ell} \right)d{\bf s} \, d{\bf t}.
\]

Since the linking number is symmetric up to sign, we may assume, without loss of generality, that $k \leq \ell$.

\begin{figure}[htbp]
	\centering
	\subfigure[ ]{\label{fig:submanifolds}\includegraphics[scale=.75]{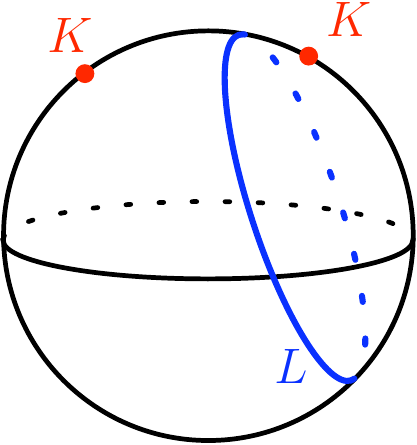}}
	\subfigure[ ]{\label{fig:cone}\includegraphics[scale=.75]{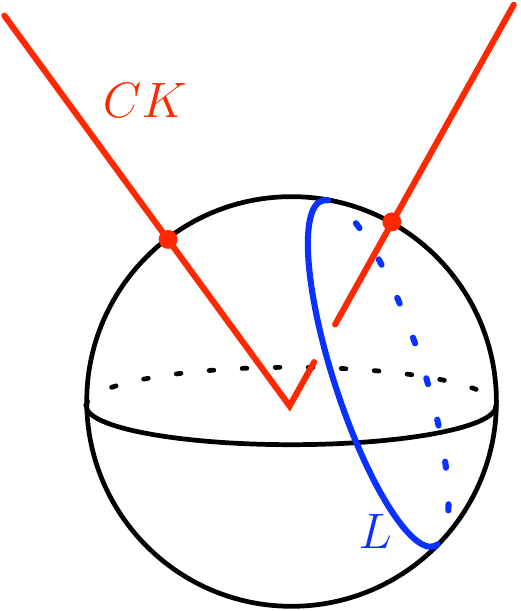}}
	\subfigure[ ]{\label{fig:cap}\includegraphics[scale=.75]{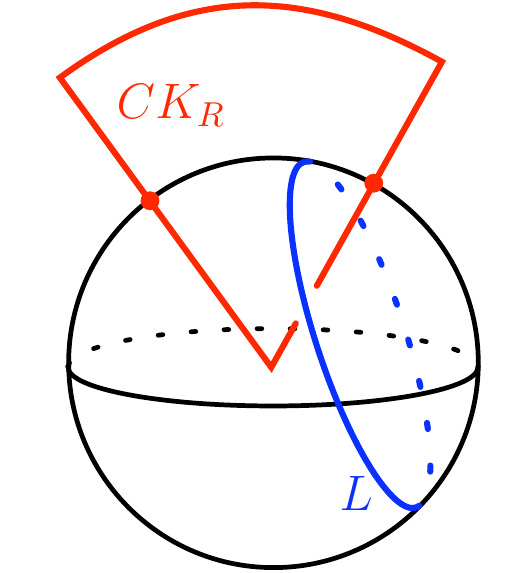}}
	\caption{ }
	\label{fig:keyfigure}
\end{figure}

The key idea is to derive an integral formula for linking numbers in $M^n$ from the one in $\mathbb{R}^{n+1}$ appearing in \thmref{thm:rN}.  This is illustrated in Figure \ref{fig:keyfigure} in the case $M^n = S^n$.  

In Figure \ref{fig:submanifolds}, we see the cycles $K$ and $L$ in $S^n$.  In Figure \ref{fig:cone}, we see the cycles $CK$ and $L$ in $\mathbb{R}^{n+1}$, where $CK := \{\tau x : x \in K, 0 \leq \tau < \infty \}$ is the (half-infinite) cone over $K$.  

We intend to take the formula for $Lk(CK,L)$ in $\mathbb{R}^{n+1}$, expressed as an integral over $CK \times L$, and partially integrate it over the rays $0 \leq \tau < \infty$ to reduce it to an integral over $K \times L$ for $Lk(K,L)$ in $S^n$.

The noncompactness of $CK$ raises questions about the meaning of $Lk(CK,L)$ and the convergence of the linking integral, which we handle by approximation as follows.  In Figure \ref{fig:cap}, we see the cycles $CK_R$ and $L$ in $\mathbb{R}^{n+1}$, where $CK_R$ consists of the truncated cone $\{\tau x : x \in K, 0 \leq \tau \leq R\}$, capped off by the scaled version $R\overline{K}$ of a chain $\overline{K}$ in $M$ bounded by $K$.  We see that $Lk_{\mathbb{R}^{n+1}}(CK_R,L) = Lk_{S^n}(K,L)$ because $CK_R$ meets a chain $\overline{L}$ bounded by $L$ exactly as $K$ does.  We then show that the integral for $Lk(CK_R,L)$ converges to the integral for $Lk(CK,L)$ as $R \to \infty$ because the contribution to the linking integral coming from the cap $R\overline{K}$ vanishes as $R \to \infty$.

\begin{lemma}\label{lem:nonclosed}
	With $K$ and $L$ as in \thmref{thm:main} and provided $n > 1$, 
	\begin{equation}\label{eqn:nonclosed}
		Lk(K, L) = \frac{(-1)^{k}}{\mathrm{vol}\, S^{n}} \int_{CK \times L} \frac{1}{|\tau x - y|^{n+1}} \, [\tau x - y, d(\tau x), d y].
	\end{equation}
\end{lemma}
\begin{proof}
	We will interpret \eqref{eqn:nonclosed} as the limit of integrals given by applying \thmref{thm:rN} to the family of closed, singular manifolds $CK_R$, $1 < R < \infty$.  By construction, $Lk(K, L) = Lk(CK_R, L)$ for all $R > 1$.  Since the singularities of $CK_R$ have measure zero, we can use \thmref{thm:rN} to compute $Lk(CK_R, L)$:
	
	\begin{align}
	\nonumber	Lk(CK_R, L) & = \frac{(-1)^{k+2}}{\mathrm{vol}\, S^{n}} \int_{CK_R \times L} \frac{1}{|\xi - y|^{n+1}} \, [\xi - y, d\xi, dy] \\
	\label{eqn:nonclosedSplit} & = \frac{(-1)^{k+2}}{\mathrm{vol}\, S^{n}} \int_{\{\tau x:x\in K, \tau \in [0,R]\}\times L} \frac{1}{|\tau x - y|^{n+1}} \, [\tau x-y, d(\tau x), dy] \\
	\nonumber & \qquad	+ \frac{(-1)^{k+2}}{\mathrm{vol}\, S^{n}}  \int_{\overline{K} \times L} \frac{1}{|R z-y|^{n+1}} \, [R z-y, d(R z), dy].
	\end{align}
	As $R \to \infty$, the first term on the right hand side of \eqref{eqn:nonclosedSplit} approaches the right hand side of \eqref{eqn:nonclosed}.  Therefore, to complete the proof we need only show that 
	\begin{equation}\label{eqn:nonclosedSecondTerm}
		\frac{(-1)^{k+2}}{\mathrm{vol}\, S^{n}} \int_{ \overline{K} \times L} \frac{1}{|R z-y|^{n+1}} \, [R z-y, d(R z), dy] \to 0
	\end{equation}
	as $R \to \infty$.  However, as $R$ gets large the numerator in the integrand grows like $R^{k+2}$, whereas the denominator grows like $R^{n+1}$.  The assumptions $k \leq \ell$ and $n > 1$ imply that $k+2 < n+1$.  Therefore, the integrand vanishes in the limit, implying \eqref{eqn:nonclosedSecondTerm} and completing the proof.
	
\end{proof}

Now we are ready to prove the linking integral formula for $M^n$.

\begin{proof}[Proof of \thmref{thm:main}]
	By \lemref{lem:nonclosed}, 
	\begin{equation}\label{eqn:mainproof1}
	Lk(K, L) = \frac{(-1)^{k}}{\mathrm{vol}\, S^{n}} \int_{CK \times L} \frac{1}{|\tau x-y|^{n+1}} \, [\tau x-y, d(\tau x), dy]. 
	\end{equation}
	Focusing first on $[\tau x-y, d(\tau x), dy]$, we see that
	\begin{align*}
		[\tau x-y, d(\tau x), dy] & = d\tau \wedge \left[\det \left(\tau x-y, x, \tau \frac{\partial x}{\partial s_1}, \ldots , \tau \frac{\partial x}{\partial s_k}, \frac{\partial y}{\partial t_1}, \ldots , \frac{\partial y}{\partial t_\ell} \right)  d{\bf s} \, d{\bf t}\right] \\
		&  = (-1)^k\tau^k d\tau \wedge [x, dx, y, dy].
	\end{align*}
	Combining this with \eqref{eqn:mainproof1}, we see that
	\[
		Lk(K, L) = \frac{1}{\mathrm{vol}\, S^{n}} \int_{K \times L} \int_{\tau = 0}^\infty \frac{\tau^k}{|\tau x -y|^{n+1}} \, d\tau \wedge [x, dx, y, dy].
	\]
	Therefore, it suffices to show that
	\begin{equation}\label{eqn:tauThetaIntegrals}
		\int_{\tau = 0}^\infty \frac{\tau^k}{|\tau x-y|^{n+1}} \, d\tau = \frac{1}{|x|^{k+1}|y|^{\ell+1}\sin^{n} \alpha} \int_{\theta = \alpha}^\pi \sin^k(\theta - \alpha) \sin^\ell\theta \, d\theta.
	\end{equation}
	\begin{figure}[htbp]
		\centering
	        \includegraphics[width=1.7in]{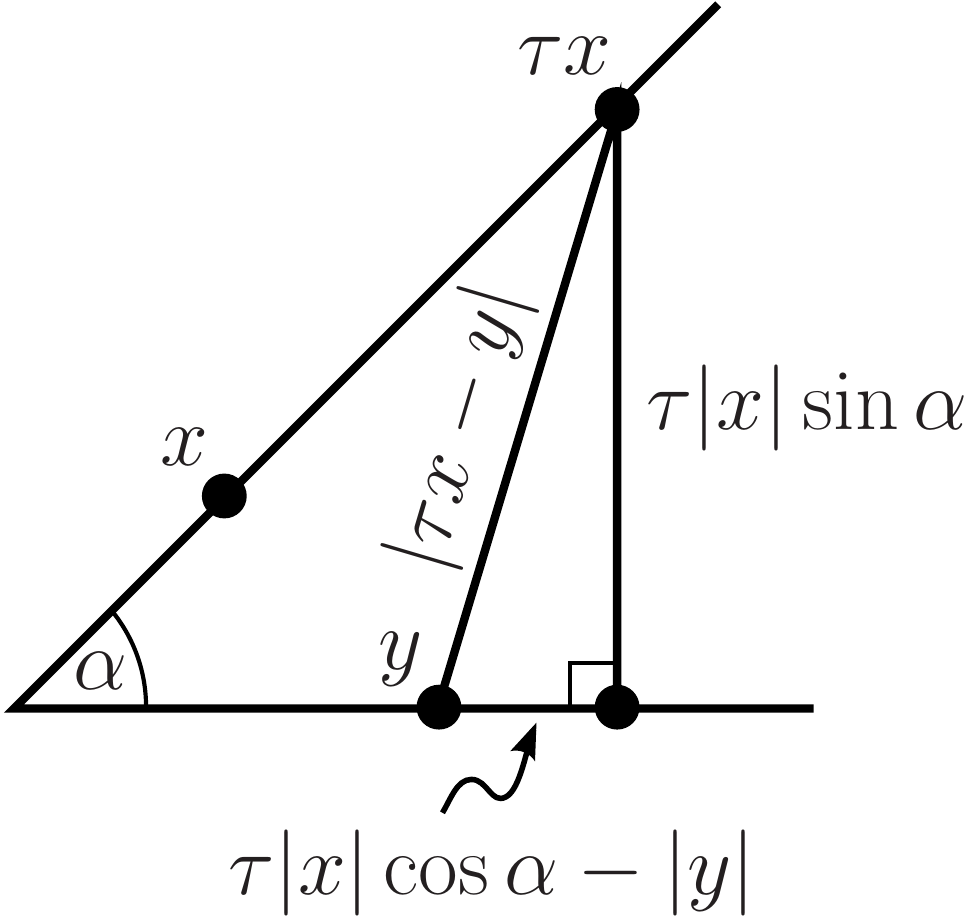}
			\caption{}
			\label{fig:triangle}
	\end{figure}

As illustrated in Figure \ref{fig:triangle}, we can re-write $|\tau x - y|^{n+1}$ as 
	\[
		|\tau x -y|^{n+1} = \left(|y|^2 + \tau^2|x|^2 - 2\tau \langle x, y \rangle \right)^{\frac{n+1}{2}} = \left(|y|^2 +  \tau^2|x|^2 - 2\tau |x||y| \cos \alpha \right)^{\frac{n+1}{2}},
	\]
where we recall that $\alpha(x,y)$ is the angle (in $\mathbb{R}^{n+1}$) formed by $x$ and $y$.

	Now make the substitution $u = \tau - \frac{|y|}{|x|} \cos \alpha$ in the left hand side of \eqref{eqn:tauThetaIntegrals} to get
	\[
		\frac{1}{|x|^{n+1}}\int_{-\frac{|y|}{|x|}\cos \alpha}^\infty \frac{\left(u + \frac{|y|}{|x|}\cos \alpha\right)^k}{\left(u^2 + \frac{|y|^2}{|x|^2}\sin^2 \alpha\right)^{\frac{n+1}{2}}} \, du.
	\]
	Substituting $-\frac{|y|}{|x|}\sin \alpha \cot \theta$ for $u$ and simplifying yields
	\begin{multline*}
		\frac{1}{|x|^{k+1}|y|^{\ell+1}\sin^{n} \alpha} \int_{\theta = \alpha}^\pi (\cos \alpha - \sin \alpha \cot \theta)^k \sin^{n-1}\theta \, d\theta \\ = \frac{1}{|x|^{k+1}|y|^{\ell+1}\sin^{n} \alpha} \int_{\theta = \alpha}^\pi \sin^k(\theta - \alpha) \sin^\ell \theta \, d\theta,
	\end{multline*}
	which is the right hand side of \eqref{eqn:tauThetaIntegrals}, completing the proof of the theorem.
\end{proof}


\section{A few examples} 
\label{sec:examples}
	To see that \thmref{thm:main} is useful in practice, we apply it to some prototypical examples.  First, we note that \thmref{thm:sn} can be used to recover the linking integral in $S^3$.  In this case, $k = \ell = 1$, so
	\begin{align*}
		\Omega_{1,1} (\alpha)& =  \int_\alpha^\pi \sin(\theta - \alpha)\sin\theta \, d\theta \\
		& = \frac{(\pi - \alpha ) \cos \alpha + \sin \alpha}{2}.
	\end{align*}
	Thus, the linking number of disjoint closed curves $K$ and $L$ in $S^3$ is given by
	\[
		Lk(K, L) = \frac{1}{4\pi^2} \int_{K \times L} \frac{(\pi - \alpha) \cos \alpha + \sin \alpha}{\sin^3 \alpha} \, [x, dx, y, dy].
	\]
	This is equivalent to the formulas obtained by DeTurck and Gluck in \cite{DG1} and Kuperberg in \cite{Kuperberg}.

	This same computation for disjoint, closed curves in $S^2 \times \mathbb{R}$ yields the linking formula
	\[
		Lk(K, L) = \frac{1}{4\pi^2} \int_{K \times L} \frac{(\pi-\alpha)\cos \alpha + \sin \alpha}{|x|^2|y|^2 \sin^3 \alpha} \, [x,dx,y,dy].
	\]
As noted in the introduction, the integrand in this linking integral formula is invariant under the action of $SO(4)$ on the ambient Euclidean space, so it is invariant under rotations in the $S^2$ factor, though not necessarily under translations in the $\mathbb{R}$ factor.  This is in contrast to the formulas of Gau\ss \ in $\mathbb{R}^3$ and DeTurck and Gluck in $S^3$ and $H^3$, whose integrands are invariant under the full group of orientation-preserving isometries of their respective manifolds.

	Now we turn to more concrete examples.  Let $K = S^k$ and $L = S^\ell$ be two great spheres contained in $S^{k+\ell + 1}$ at constant geodesic distance $\pi/2$.  Up to ambient isometry, we may as well take $K$ as the unit sphere in the first $k+1$ coordinates in $\mathbb{R}^{k+\ell+2}$ and $L$ as the unit sphere in the last $\ell + 1$ coordinates.  Then, by construction, $Lk(K, L) = 1$; we want to see that \thmref{thm:sn} gives the same result.
	
	To simplify the computation, we recall the following fact relating integrals of sines and cosines to the volumes of spheres:
	\begin{fact}\label{fact:sinecosinevolspheres}
		\[
			\mathrm{vol}\, S^{k+\ell+1} = \mathrm{vol}\, S^k \ \mathrm{vol}\, S^\ell\int_0^{\pi/2} \cos^k \theta \sin^\ell \theta \, d\theta .
		\]
	\end{fact}
	Fact \ref{fact:sinecosinevolspheres} follows from viewing $S^{k+\ell+1}$ as the spherical join of orthogonal great spheres $S^k$ and $S^\ell$ and computing the volume in the natural coordinates that result from this identification.

	In our example, since $K$ and $L$ are at constant geodesic distance $\pi/2$, the above fact implies that $\Omega_{k, \ell}$ has the form
	\begin{align*}
		\Omega_{k, \ell} (\alpha)& = \int_{\pi/2}^\pi \sin^k (\theta - \pi/2) \sin^\ell \theta \, d\theta \\
		& = \int_0^{\pi/2} \cos^k \theta \sin^\ell \theta \, d\theta \\
		& = \frac{\mathrm{vol}\, S^{k+\ell+1}}{\mathrm{vol}\, S^k \ \mathrm{vol}\, S^\ell}.
	\end{align*}
	Therefore, 
	\begin{align*}
		Lk(K, L) & = \frac{1}{\mathrm{vol}\, S^{k+\ell+1}} \int_{K \times L} \frac{1}{\sin^{k+\ell+1} (\pi/2)}  \frac{\mathrm{vol}\, S^{k+\ell+1}}{\mathrm{vol}\, S^k \ \mathrm{vol}\, S^\ell} \, [x,dx,y,dy] \\
		& = \frac{1}{\mathrm{vol}\, S^k \ \mathrm{vol}\, S^\ell} \int_{K\times L} [x,dx,y,dy] \\
		& = \frac{1}{\mathrm{vol}\, S^k \ \mathrm{vol}\, S^\ell} \mathrm{vol}\, S^k \ \mathrm{vol}\, S^\ell \\
		& = 1.
	\end{align*}
	The second-to-last equality follows because, since $K$ is contained in the  $\mathbb{R}^{k+1}$ factor and $L$ is contained in the $\mathbb{R}^{\ell+1}$ factor of $\mathbb{R}^{k+\ell+2} = \mathbb{R}^{k+1} \times \mathbb{R}^{\ell+1}$, $[x,dx,y,dy] = [x,dx] \wedge [y, dy] = d\mr{Vol}_{S^k} \wedge d\mr{Vol}_{S^\ell}$.
	
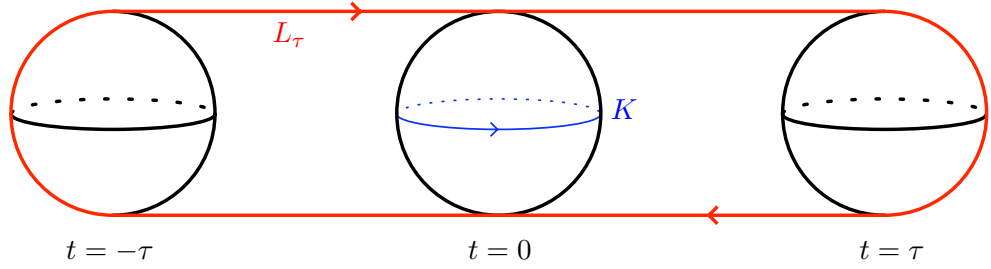
\begin{figure}[htbp]
		\centering
		\begin{picture}(370,92)
			\put(0,12){\includegraphics[width=370px]{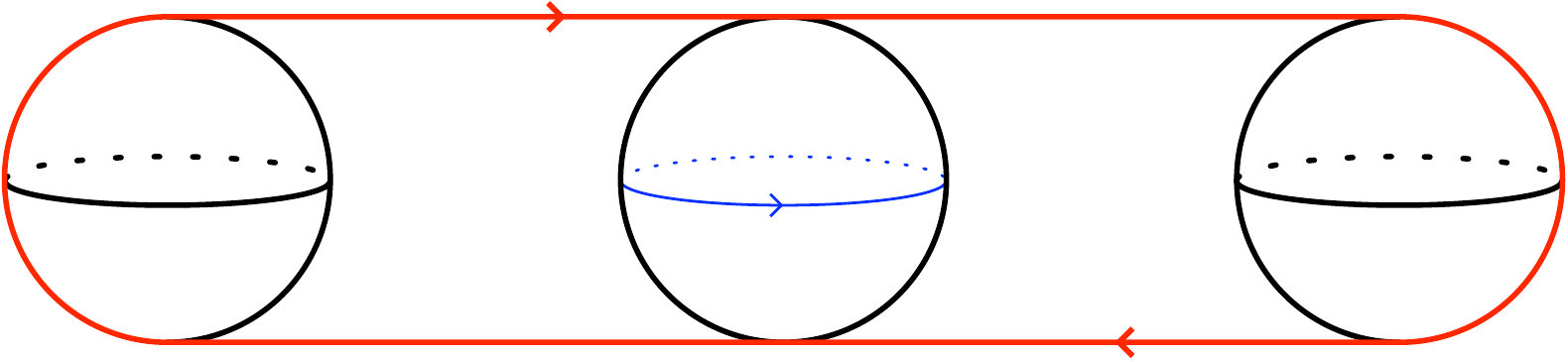}}
			\put(228, 53){\textcolor{blue}{$K$}}
			\put(100, 82){\textcolor{red}{$L_\tau$}}
			\put(22, 0){$t = -\tau$}
			\put(174, 0){$t = 0$}
			\put(322, 0){$t = \tau$}
		\end{picture}
		\caption{$K$ and $L_\tau$ in $S^2 \times \mathbb{R}$}
		\label{fig:s2xrloop}
	\end{figure}

	Finally, we use \thmref{thm:main} to compute the simplest example in $S^2 \times \mathbb{R}$.  Thinking of the $\mathbb{R}$ factor as time, let $K$ be the instantaneous equator on the $S^2$ factor at time 0 and let $L$ be the union of the eternal north pole traversed positively in time with the eternal south pole traversed negatively in time.  We can think of $L$ as the limit when $\tau \to \infty$ of the curves $L_\tau$ shown in Figure \ref{fig:s2xrloop}. 
	
	  With this choice of $K$ and $L$, $\alpha$ is again equal to $\pi/2$ for all $x\in K$ and $y\in L$.  Therefore, 
	\[
		\Omega_{1,1}(\alpha) = \frac{\mathrm{vol}\, S^3}{\mathrm{vol}\, S^1 \ \mathrm{vol}\, S^1} = \frac{\mathrm{vol}\, S^3}{4\pi^2}.
	\]
	Thus,
	\begin{align*}
		Lk(K, L) & = \frac{1}{\mathrm{vol}\, S^3} \int_{K \times L} \frac{\mathrm{vol}\, S^3}{|y|^2 4\pi^2 } \, [x,dx,y,dy] \\
		& = 2 \cdot \frac{1}{4\pi^2}  \int_{-\infty}^\infty \int_0^{2\pi}\frac{1}{1+t^2} \det \left(\begin{array}{cccc} \cos s & \sin s & 0 & 0 \\  -\sin s & \cos s & 0 & 0 \\ 0 & 0 & 1 & t \\ 0 & 0 & 0 & 1 \end{array} \right) ds\, dt \\
		& = \frac{1}{2\pi^2} \int_{-\infty}^\infty \int_{0}^{2\pi} \frac{1}{1+t^2}\, ds\, dt \\
		& = 1.
	\end{align*}
	In the second equality we use the symmetry of $L$ to eliminate the integral over the south pole component by doubling the contribution of the north pole component. 

\end{document}